\documentclass[12pt]{article}
\usepackage{amssymb}
\newtheorem{theorem}{Theorem}
\newtheorem{corollary}[theorem]{Corollary}
\newtheorem{definition}{Definition}
\newtheorem{lemma}[theorem]{Lemma}
\newtheorem{proposition}[theorem]{Proposition}

\newenvironment{proof}{\begin{trivlist}\item[]{\it
Proof.}}{\hfill$\square$\end{trivlist}}

\newenvironment{proofofprop}[1]{\noindent{\it Proof of Proposition
#1}}{\hfill$\square$\\\mbox{}}

\def\span{{\mathrm{Span}}}
\def\char{{\mathrm{char}}}
\def\sltwo{{\mathrm{SL}}(2,{\mathbb{C}})}
\def\mc{{\mathbb{C}}}
\def\onepsg{1{\mathrm{-PS}}}
\def\limzr{{\mathrm{lim}}_{z\to 0}\rho(z)}
\def\pol{{\mathrm{Pol}}}

\begin{document}
\title{Typical separating invariants}
\author{M. Domokos
\\ 
\\ 
{\small R\'enyi Institute of Mathematics, Hungarian Academy of 
Sciences,} 
\\ {\small P.O. Box 127, 1364 Budapest, Hungary,} 
{\small E-mail: domokos@renyi.hu } }
\date{}
\maketitle 
\begin{abstract} 
It is shown that a trivial version of polarization is sufficient to produce separating systems of polynomial invariants: if two points in the direct sum of the $G$--modules  
$W$ and $m$ copies of $V$ can be separated by polynomial invariants, then 
they can be separated by invariants depending only on $\leq 2\dim(V)$ variables  
of type $V$; when $G$ is reductive, invariants depending only on 
$\leq \dim(V)+1$ variables suffice. Similar result is valid for rational invariants. Explicit bounds on the number of type $V$ variables in a typical system of separating invariants are given for the binary 
polyhedral groups, and this is applied to the invariant theory of binary forms. 
\end{abstract}


\section{Introduction}\label{sec:intro}

Let $G$ be a group, $k$ an arbitrary base field, $V,W$ finite dimensional 
$G$--modules over $k$. 
Write $V^m$ for the $m$--fold direct sum $V\oplus\cdots\oplus V$. 
The coordinate ring $k[V]$ of $V$ 
(the symmetric tensor algebra of the dual $G$--module $V^*$) 
contains the subalgebra $k[V]^G$ of 
polynomial invariants. We say that $v,v'\in V$ {\it can be separated 
(by polynomial invariants)} if there exists an $f\in k[V]^G$ with 
$f(v)\neq f(v')$. 
Following \cite{dk}, Definition 2.3.8, we say that 
$S\subseteq k[V]^G$ is a {\it separating set} if whenever 
$v,v'\in V$ can be separated, then there exists an $f\in S$ with 
$f(v)\neq f(v')$. 

We shall be interested in separating sets in 
$k[W\oplus V^m]^G$, where $m$ is varying in $\mathbb{N}$. 
Given $f\in k[W\oplus V^d]$, $m\geq d\in\mathbb{N}$, 
and $1\leq i_1<\cdots<i_d\leq m$, denote by 
$f^{(i_1,\ldots,i_d)}$ the function 
$f\circ\pi_{(i_1,\ldots,i_d)}:W\oplus V^m\to k$, where 
$\pi_{(i_1,\ldots,i_d)}$ is the projection  
mapping $(w,v_1,\ldots,v_m)\mapsto (w,v_{i_1},\ldots,v_{i_d})$. 

\begin{definition}\label{def:cstsi} 
We say that $S\subseteq k[W\oplus V^d]^G$ is a 
{\it complete system of typical separating invariants in type $V$ 
variables relative to $W$}, if for all $m\in \mathbb{N}$, 
$m\geq d$, we have that 
\[S^{(m)}=\{f^{(i_1,\ldots,i_d)}\mid f\in S,\quad 
1\leq i_1<\cdots <i_d\leq m\}\] 
is a separating set in $k[W\oplus V^m]^G$. 
(Note that in this case for $q<d$, substituting zeros into the type $V$ 
variables of index $q+1,q+2,\ldots,d$ in the elements of $S$ we obtain 
a separating set in $k[W\oplus V^q]^G$.) 
\end{definition}

This concept can be considered as a naive version of polarization, as 
the elements of $S^{(m)}$ are very special polarizations of the elements of 
$S$ (see Section~\ref{sec:conclude} for comments on polarization). 
Whereas the general polarization process is a basic tool (at least in characteristic zero) to produce generating invariants, the analogues of the above concept for generating systems seems to be useless. 
However, it turns out to be meaningful when we deal with separating sets. 

\begin{definition}\label{def:sigma} 
Denote by $\sigma(G,W,V)$ the minimal non-negative integer $d$ such that 
$k[W\oplus V^d]^G$ contains a complete system of typical separating invariants in type $V$ variables relative to $W$. 
\end{definition} 

In Section~\ref{sec:general} we verify the general bound 
$\sigma(G,W,V)\leq 2\dim(V)$. 
For a reductive group $G$ this is improved in Section~\ref{sec:reductive} 
to $\sigma(G,W,V)\leq 1+\dim(V)$; the latter bound is sharp. 
We prove in Section~\ref{sec:finite} that the supremum of the numbers $\sigma(G,W,V)$ for a fixed finite group $G$ can be interpreted as a purely 
group theoretic invariant that we call the Helly dimension of $G$. 
In Section~\ref{sec:binary} we turn to the invariant theory of binary forms: we investigate the numbers $\sigma(G,W,V)$ when 
$G$ is the special linear group $\sltwo$ and $V$ is an irreducible 
$\sltwo$--module. It turns out that there is a universal constant upper bound 
for $\sigma(G,W,V)$, independent of $W$ and $V$. 
In Section~\ref{sec:rational} we apply the idea of Definition~\ref{def:cstsi} 
for rational invariants (where the general polarization method can not be applied). 


\section{General bound for the number of variables}\label{sec:general} 

Denote by $Gr(V\oplus V)$ the set of linear subspaces in $V\oplus V$. 

\begin{lemma}\label{lemma:subspace} 
Whether $v=(w,v_1,\ldots,v_m)$ and 
$v'=(w',v_1',\ldots, v_m')$ can be separated or not depends 
only on the triple 
$(w,w',\span_k\{v_1\oplus v_1',\ldots,v_m\oplus v_m'\})
\in W\times W\times Gr(V\oplus V)$. 
\end{lemma} 

\begin{proof} 
Take $v=(w, v_1,\ldots,v_m)\in W\oplus V^m$, 
$v'=(w',v_1',\ldots,v_m')\in W\oplus V^m$, 
$u=(w,u_1,\ldots,u_d)\in W\oplus V^d$, 
$u'=(w',u_1',\ldots,u_d')\in W\oplus V^d$ 
with 
\[\span_k(v_1\oplus v_1',\ldots,v_m\oplus v_m'\}
=\span_k(u_1\oplus u_1',\ldots,u_d\oplus u_d'\}.\] 
In particular, there are coefficients $\alpha_{ij}\in k$ 
such that 
\[v_i\oplus v_i'=\sum_{j=1}^d \alpha_{ij}(u_j\oplus u_j')\quad 
\mbox{for}\quad i=1,\ldots,m.\] 
Given an invariant $f\in k[W\oplus V^m]^G$, define 
$h\in k[W\oplus V^d]$ by 
\[h(x,y_1,\ldots,y_d)=
f(x,\sum_{j=1}^d\alpha_{1j}y_j,\ldots,\sum_{j=1}^d\alpha_{mj}y_j),\] 
where $x$ (resp. $y_j$) stand for vector variables of type $W$ 
(resp. $V$). 
Obviously, $h$ is $G$--invariant. By construction we have 
$h(u)=f(v)$ and $h(u')=f(v')$. So if $f(v)\neq f(v')$, then 
$h(u)\neq h(u')$. Consequently, if $v$ and $v'$ can be separated, 
then $u$ and $u'$ can be separated. 
The reverse implication is shown in the same way. 
\end{proof}

\begin{theorem}\label{thm:2n}
Let $S$ be a separating set in $k[W\oplus V^{2n}]^G$, 
where $n=\dim(V)$. 
Then $S$ is a complete system of typical separating invariants in type $V$ variables relative to $W$. 
That is, we have $\sigma(G,W,V)\leq 2\dim(V)$. 
\end{theorem}

\begin{proof} 
Take a pair of points $v=(w,v_1,\ldots,v_m)$, 
$v'=(w',v_1',\ldots,v_m')$ in $W\oplus V^m$ where $m\geq 2n$. 
Since $\dim(\span_k\{v_1\oplus v_1',\ldots,v_m\oplus v_m'\})$ 
is at most $2n$, there are indices 
$1\leq i_1<\cdots <i_{2n}\leq m$ with 
\[\span_k\{v_{i_1}\oplus v_{i_1}',\ldots,v_{i_{2n}}\oplus v_{i_{2n}}'\}
=\span_k\{v_1\oplus v_1',\ldots,v_m\oplus v_m'\}.\] 
If $v$ and $v'$ can be separated, then 
$(w,v_{i_1},\ldots,v_{i_{2n}})$ and 
$(w',v_{i_1}',\ldots,v_{i_{2n}}')$ 
can be separated by Lemma~\ref{lemma:subspace}, so there exists an $f\in S$ 
separating them. This means that 
$f^{(i_1,\ldots,i_{2n})}\in k[W\oplus V^m]^G$ separates $v$ and $v'$. 
\end{proof}


\section{Reductive groups}\label{sec:reductive} 

Throughout this section we assume that $k$ is algebraically closed. 
It is well known that if 
$G$ is reductive (acting algebraically on $V$), then two points in $V$ can be separated if and only if the Zariski closures of their $G$--orbits are disjoint, see for example Corollary 3.5.2 in \cite{n}. 
Recall also that the Zariski closure of any $G$--orbit contains a unique closed orbit, hence a set of invariants is separating if it separates closed orbits. 

\begin{lemma}\label{lemma:closed}
Suppose that $v=(w,v_1,\ldots,v_m)\in W\oplus V^m$, 
and for some $q\leq m$ we have that 
$\span_k\{v_1,\ldots,v_q\}=\span_k\{v_1,\ldots,v_m\}$. 
Then the $G$--orbit of 
$v$ is closed in 
$W\oplus V^m$ if and only if 
the $G$--orbit of $(w,v_1,\ldots,v_q)$ 
is closed in $W\oplus V^q$. 
\end{lemma} 

\begin{proof} 
There are coefficients $\beta_{q+i,j}\in k$ with 
$v_{q+i}=\sum_{j=1}^q\beta_{q+i,j}v_j$ holding for $i=1,\ldots,m-q$. 
The map 
\[\varphi:W\oplus V^q\to W\oplus V^m\] 
\[(x,y_1,\ldots,y_q)\mapsto 
(x,y_1,\ldots,y_q,\sum_{j=1}^q\beta_{q+1,j}y_j, 
\ldots,\sum_{j=1}^q\beta_{m,j}y_j)\] 
is an injective morphism of affine $G$--varieties, its image 
is a $G$--stable linear subspace $L$ of $W\oplus V^m$. 
So $\varphi$ induces an isomorphism of the affine $G$-varieties 
$W\oplus V^q$ and $L$, moreover, $\varphi$ maps 
the $G$--orbit of $(w,v_1,\ldots,v_q)$ in $W\oplus V^q$ onto  
the $G$--orbit of $v$ in $W\oplus V^m$. 
Hence the $G$--orbit of $(w,v_1,\ldots,v_q)$ 
is closed in $W\oplus V^q$ if and only if the $G$--orbit of $v$ is closed in $L$. Trivially, a subset of $L$ is closed in $L$ 
if and only if it is closed in $W\oplus V^m$. 
\end{proof}

\begin{theorem}\label{thm:n+1} 
Assume that $G$ is reductive, and 
let $S$ be a separating set in $k[W\oplus V^{n+1}]^G$, where 
$n=\dim(V)$. 
Then $S$ is a complete system of typical separating invariants 
in type $V$ variables relative to $W$. 
In other words, $\sigma(G,W,V)\leq \dim(V)+1$ if $G$ is reductive. 
\end{theorem} 

\begin{proof} 
Take $v=(w,v_1,\ldots,v_m)\in W\oplus V^m$ and 
$v'=(w',v_1',\ldots,v_m')\in W\oplus V^m$ where 
$m\geq n+1$, and assume that the $G$--orbits of $v$ and $v'$ are different and closed with respect to the Zariski topology. We need to show that after a possible rearrangement of the type $V$ variables, $v$ and $v'$ can be separated by an element of $S$ (identifying $f\in S$ with $f^{(1,\ldots,n+1)}$). 

Without loss of generality we may assume that 
\[\dim(\span_k\{v_1,\ldots,v_m\})\geq 
\dim(\span_k\{v_1',\ldots,v_m'\}),\] 
and that 
$v_1,\ldots,v_q$ is a basis of $\span_k\{v_1,\ldots,v_m\}$ for some 
$q\leq n$. 

{\it Case I. $\span_k\{v_1',\ldots,v_q'\}\neq \span_k\{v_1',\ldots,v_m'\}$} 

Our assumptions imply that $v_1',\ldots,v_q'$ are linearly dependent, 
so $\lambda_1 v_1'+\cdots+\lambda_qv_q'=0$ for some $\lambda_i\in k$, not all zero. 
It follows that the closure of the $G$--orbit of 
$(w',v_1',\ldots,v_q')$ is contained in the 
proper linear subspace $L$ of $W\oplus V^q$ defined by the system of linear equations 
$\sum_{i=1}^q\lambda_iy_i=0$. On the other hand, the $G$--orbit of 
$(w,v_1,\ldots,v_q)$ is not contained in $L$, and is closed in 
$W\oplus V^q$ by Lemma~\ref{lemma:closed}. 
Consequently, $(w,v_1,\ldots,v_q)$ and 
$(w',v_1',\ldots,v_q')$ can be separated by an element 
of $k[W\oplus V^q]^G$, implying that $v$ and $v'$ can be separated by an element of $S$. 

{\it Case II. $\span_k\{v_1',\ldots,v_q'\}= \span_k\{v_1',\ldots,v_m'\}$} 

By Lemma~\ref{lemma:closed}, the $G$--orbits of 
$(w,v_1,\ldots,v_q)$ and $(w',v_1',\ldots,v_q')$ 
are both closed in $W\oplus V^q$. If these $G$--orbits are different, then 
they can be separated by an element of $k[W\oplus V^q]^G$, whence 
they can be separated by an element of $S$. 
Otherwise replace $v'$ by $g\cdot v'$ with a $g\in G$ satisfying 
$g\cdot v_i'=v_i$ for $i=1,\ldots,q$ and $g\cdot w'=w$. 
In other words, we may assume from now on that $w=w'$, 
$v_1=v_1',\ldots,v_q=v_q'$, and 
$v_1,\ldots,v_q$ is a basis of 
$\span_k\{v_1,\ldots,v_m\}=\span_k\{v_1',\ldots,v_m'\}$. 
There are unique coefficients $\gamma_{ij},\gamma_{ij}'\in k$ 
such that 
\[\sum_{j=1}^q\gamma_{ij}v_j=v_{q+i}, \quad 
\sum_{j=1}^q\gamma_{ij}'v_j'=v_{q+i}'\quad (i=1,\ldots,m-q).\] 
Since $v$ and $v'$ are different, $v_{q+i}\neq v_{q+i}'$ for some $i$, 
say $v_{q+1}\neq v_{q+1}'$, implying that 
$\gamma_{1,j}\neq\gamma_{1,j}'$ for some $j\in\{1,\ldots,q\}$. 
It follows that the $G$--orbits of 
$(w,v_1,\ldots,v_{q+1})$ and 
$(w',v_1',\ldots,v_{q+1}')$ are different in $W\oplus V^{q+1}$ 
(since the action of $G$ preserves the linear dependency relations among the type $V$ variables). 
These orbits are closed in $W\oplus V^{q+1}$ by Lemma~\ref{lemma:closed}, 
hence they can be separated by an invariant in $k[W\oplus V^{q+1}]^G$. 
Now $q\leq n$, so this shows that $v$ and $v'$ (after a possible rearrangement 
of the type $V$ variables) 
can be separated by an element of $S$.  
\end{proof}

{\bf Remark.} Similar arguments imply that the nullcone is defined by invariants depending only on $\dim(V)$ variables of type $V$.

The following example shows that the bound $\dim(V)+1$ on the 
number of type $V$ variables in Theorem~\ref{thm:n+1} is sharp. 

{\bf Example.} Let $G$ be a torus $(k^{\times})^{n-1}$ of rank $n-1$ acting on 
$V=k^n$ by the rule 
\[(\alpha_1,\ldots,\alpha_{n-1})\cdot 
\left[\begin{array}{c}x_1 \\x_2 \\x_3 \\\vdots \\x_n\end{array}\right]
=
\left[\begin{array}{c}(\alpha_1\cdots\alpha_{n-1})x_1 \\\alpha_1^{-2}x_2 \\\alpha_2^{-2}x_3 \\\vdots \\\alpha_{n-1}^{-2}x_n\end{array}\right].\] 
Denote by $e_1,\ldots,e_n$ the standard basis vectors in $k^n$, and set 
$v=(-e_1,e_1,e_2,e_3,\ldots,e_n)$ 
and 
$v'=(e_1,e_1,e_2,e_3,\ldots,e_n)$. 
Write $x(i)_j$ for the coordinate function on $V^{n+1}$ mapping an 
$(n+1)$--tuple of vectors to the $j$th coordinate of the $i$th component. 
Then 
\[x(1)_1x(2)_1x(3)_2x(4)_3x(5)_4\cdots x(n+1)_n\] 
is an invariant in 
$k[V^{n+1}]^G$ separating $v$ and $v'$ (provided that the characteristic of 
$k$ is different from $2$). On the other hand, we show that $v$ and $v'$ can not be separated by invariants depending only on $\leq n$ vector variables. 
Define the multidegree of a monomial in $k[V^{n+1}]$ to be 
$(d_1,\ldots,d_n)$ if its total degree in the variables 
$x(1)_i,\ldots,x(n+1)_i$ is $d_i$ for $i=1,\ldots,n$. It is easy to see that 
an element of $k[V^{n+1}]$ is $G$--invariant if and only if all of its non-zero monomials have multidegree of the form $(2d,d,d,\ldots,d)$ with $d\in\mathbb{N}$. 
Now let $f$ be such a monomial. Then either $f(v)=0=f(v')$, or 
\[f=x(1)_1^bx(2)_1^{2d-b}(x(3)_2x(4)_3x(5)_4\cdots x(n+1)_n)^d 
\quad (b\in\{0,1,\ldots,2d\}).\] 
In the latter case, if $f$ involves only $n$ vector variables, then 
$b=0$ or $b=2d$, implying that 
$f(-e_1,u)=f(e_1,u)$ for all $u\in V^n$. 
If the characteristic of the base field is $2$, then we modify the action to 
\[(\alpha_1,\ldots,\alpha_{n-1})\cdot 
\left[\begin{array}{c}x_1 \\x_2 \\x_3 \\\vdots \\x_n\end{array}\right]
=
\left[\begin{array}{c}(\alpha_1\cdots\alpha_{n-1})x_1 \\\alpha_1^{-3}x_2 \\\alpha_2^{-3}x_3 \\\vdots \\\alpha_{n-1}^{-3}x_n\end{array}\right],\] 
and show in the same way as above that $(\omega e_1,e_1,e_2,\ldots,e_n)$ (where $\omega$ is a primitive third root of $1$) 
and $(e_1,e_1,e_2,\ldots,e_n)$ can not be separated by invariants depending on $\leq n$ vector variables. 


\section{Finite groups}\label{sec:finite} 

Throughout this section we assume that the group $G$ is finite. 
Then all orbits are closed for any linear action of $G$, and two points in a 
$G$--module can be separated by polynomial invariants if and only if their 
orbits are different. 
This implies upper bounds for $\sigma(G,W,V)$ given 
in terms of the subgroup lattice of $G$. 

Define the {\it Helly dimension} $\kappa(G)$ of $G$ as the minimal 
natural number $n$ with the following property: if we are given $m\geq n$ 
left cosets in $G$ such that any $n$ of them have a non-empty 
intersection, then the $m$ cosets have a common element. 
We say that the subgroups $G_1,\ldots,G_m$ are {\it intersection independent} 
if no $G_i$ contains $\cap_{j\neq i}G_j$, the intersection of the others. 
Note that in this case $G_1,\ldots,G_m$ are necessarily distinct, 
and if $m>1$, then the $G_i$ are proper subgroups. Denote by $\mu(G)$ the maximal size of an intersection independent set of subgroups of $G$, and write $\lambda(G)$ 
for the maximal length of a chain of proper subgroups in $G$. 

{\bf Remark.} The name "Helly dimension" is motivated by Helly's Theorem (cf. \cite{h}). Numbers related to $\kappa(G)$ have already been introduced and studied for certain finite groups, see for example \cite{cc}. 

\begin{lemma}\label{lemma:inequalities} 
We have the inequalities 
$\kappa(G)\leq \mu(G)+1$ and $\mu(G)\leq \lambda(G)$.  
\end{lemma} 

\begin{proof} Take an $n\geq \mu(G)+1$ and 
$n+1$ left cosets $g_0G_0,g_1G_1,\ldots,g_nG_n$ in $G$ such that any 
$n$ of them have a non-empty intersection. 
Choose an element $g$ from $\cap_{i=1}^ng_iG_i$. 
Since $G_1,\ldots,G_n$ are not intersection independent, one of them, say $G_n$ 
contains $\cap_{i=1}^{n-1}G_i.$ By assumption there is an element $h$ 
in $\cap_{i=0}^{n-1}g_iG_i=g_0G_0\cap g(\cap_{i=1}^{n-1}G_i)$, 
so $h$ is contained also in $gG_n=g_nG_n$, showing that 
the intersection of all the cosets $g_iG_i$ $(i=0,\ldots,n)$ is non-empty. 
This proves the first inequality. 

For the second inequality observe that if $G_1,\ldots,G_n$ is intersection independent, then 
$G_1\supset G_1\cap G_2\supset G_1\cap G_2\cap G_3\supset\cdots\supset 
\cap_{i=1}^n G_i$ 
is a strictly descending chain of subgroups, and if $n>1$, then 
$G_1$ is proper. 
\end{proof}

The Helly dimension is a very natural characteristic of a group. 
Given an action of $G$ on some set $X$, $d\leq m\in \mathbb{N}$, 
consider the diagonal action on $X^m=X\times\cdots\times X$, 
and ask whether two points $x=(x_1,\ldots,x_m)$ and 
$x'=(x_1',\ldots,x_m')$ belong to the same orbit, provided that 
$(x_{i_1},\ldots,x_{i_d})$ and $(x_{i_1}',\ldots,x_{i_d}')$ belong to the same orbit in $X^d$ for all $1\leq i_1<\cdots<i_d\leq m$. 
The minimal $d$ such that this holds for all actions of $G$ and all $m\geq d$ 
is nothing but $\kappa(G)$. 
The following variant of this observation explains our interest in 
$\kappa(G)$. 

\begin{proposition}\label{prop:helly} 
The Helly dimension $\kappa(G)$ is the maximal value of 
$\sigma(G,W,V)$, as $W,V$ range over all finite dimensional $G$--modules. 
\end{proposition} 

\begin{proof} Recall that two points in a $G$--module can not be separated 
by polynomial invariants if and only if their orbit is the same. 

First we prove the inequality 
$\sigma(G,W,V)\leq\kappa(G)$ for any $G$--modules $W,V$. 
Assume that $m\geq d=\kappa(G)$, 
$v=(w,v_1,\ldots,v_m)\in W\oplus V^m$, 
$v'=(w',v_1',\ldots,v_m')\in W\oplus V^m$, 
and suppose that for all $1\leq i_1<\cdots<i_d\leq m$ we have that 
$\pi_{(i_1,\ldots,i_d)}(v)$ and 
$\pi_{(i_1,\ldots,i_d)}(v')$ 
belong to the same $G$--orbit in $W\oplus V^d$. We need to show that 
$v$ and $v'$ belong to the same $G$--orbit. 
Replacing $v$ by an appropriate element in its orbit we may assume 
that $w=w'$. The set of elements $g$ in the stabilizer of $w$ 
with the property $gv_i=v_i'$ is a left coset $gG_i$ in $G$, for 
$i=1,\ldots,m$. 
Our assumption is that any $d$ of these cosets have a non-empty intersection. 
Since $d=\kappa(G)$, it follows that there is a common element $g$ of these cosets, and obviously $g\cdot v=v'$. 

To show the reverse inequality, take a natural number $d<\kappa(G)$. 
Then there is an $m>d$ and left cosets $g_1G_1,\ldots,g_mG_m$ such that 
any $d$ of them have a non-empty intersection, but the $m$ cosets have no common element. It is well known that each subgroup of a finite group can be realized 
as a stabiliser of a vector in an appropriate finite dimensional representation (consider the representation induced by the trivial representation of the subgroup). Taking direct sums one can construct a finite 
dimensional $G$--module $V$ and vectors $v_1,\ldots,v_m\in V$ such that the 
stabilizer of $v_i$ is $G_i$ for $i=1,\ldots,m$. 
Set $v=(v_1,\ldots,v_m)\in V^m$ and $v'=(g_1v_1,\ldots,g_mv_m)\in V^m$. 
Our assumptions imply that $\pi_{(i_1,\ldots,i_d)}(v)$ and 
$\pi_{(i_1,\ldots,i_d)}(v')$ belong to the same $G$--orbit 
for all $1\leq i_1<\cdots<i_d\leq m$, but $v$ and $v'$ belong to 
different orbits. This shows that $\sigma(G,0,V)>d$. 
Since this holds for all $d<\kappa(G)$, we conclude that the supremum 
of the numbers $\sigma(G,W,V)$ as $W$ and $V$ vary is not smaller than 
$\kappa(G)$. 
\end{proof}

In the rest of this section we work out an explicit bound for the Helly dimension of the finite subgroups of $\sltwo$ (the binary polyhedral groups); this will be used in Section~\ref{sec:binary}.   
Recall that up to isomorphism the finite subgroups of $SO_3$ (the group of rotations of the $3$--dimensional real Euclidean space) are the following: 
the cyclic group $C_n$ of order $n$ (the group of rotations of a regular $n$--angle pyramid), the dihedral group $D_n$ of order $2n$ (the group of rotations stabilizing a regular $n$--gon), 
the alternating group $A_4$ (the group of rotations of a regular tetrahedron), 
the symmetric group $S_4$ (the group of rotations of a cube), and the alternating group $A_5$ (the group of rotations of a regular icosahedron). 
Denote by $\widetilde D_n$, $\widetilde A_4$, $\widetilde S_4$, 
$\widetilde A_5$ the preimages of these groups under the natural 
double covering $SU_2\to SO_3$. Now 
up to isomorphism the following is a complete list of finite subgroups of 
$\sltwo$ (see for example Chapter 0.13 in \cite{pv}): 
\[C_n,\quad \widetilde D_n \quad(n=1,2,\ldots), \quad \widetilde A_4, \quad 
\widetilde S_4, \quad \widetilde A_5.\] 
 
\begin{proposition}\label{prop:cyclic} 
We have $\kappa(C_n)=2$ for any $n>1$. 
\end{proposition} 

\begin{proof} 
When $n>1$, consider the cosets of two different elements with respect to 
the trivial subgroup $\{1\}$; this example shows that $\kappa(G)\geq 2$. 
To show the reverse inequality, denote by $g$ a generator of $C_n$. 
Any coset in $C_n$ is of the form $G_{q,r}=\{g^x\mid x\equiv r\mbox{ mod }q\}$ 
for some positive divisor $q$ of $n$ and $r\in\{0,1,\ldots,q-1\}$. 
Take $m\geq 2$ and $m$ cosets $G_{q_1,r_1},\ldots,G_{q_m,r_m}$ such that 
any pair has a non-empty intersection. Consider the system of congruences 
\[x\equiv r_1\mbox{ mod }q_1;\ 
x\equiv r_2 \mbox{ mod }q_2;\ \ldots;\  
x\equiv r_m\mbox{ mod }q_m. \]
By assumption, any pair of these congruences has a simultaneous solution.  
Therefore the system has a solution by the Chinese Remainder Theorem. 
This means that the intersection of the $m$ given cosets is non-empty. 
\end{proof} 

{\bf Remark.} For comparison we mention that $\mu(C_n)$ can be 
arbitrary large: for example, if $n$ is the product of $r$ distinct primes 
$p_1,\ldots,p_r$, then the subgroups of order 
$\frac{n}{p_1},\frac{n}{p_2},\ldots,\frac{n}{p_r}$ 
are intersection independent. 

\begin{proposition}\label{prop:dihedral} 
We have $\kappa(\widetilde D_n)\leq 4$ for any $n\geq 2$. 
\end{proposition}
\begin{proof} 
In terms of generators and relations, 
$\widetilde D_n=\langle a,b\mid b^4=1,\quad a^n=b^2,\quad 
bab^{-1}=a^{-1}\rangle$. 
For a positive divisor $q$ of $2n$ and $r\in\{0,1,\ldots,q-1\}$, 
set 
$A_{q,r}=\{a^x\mid x\equiv r\mbox{ mod }q\}$ and 
$B_{q,r}=\{ba^x\mid x\equiv r\mbox{ mod }q\}$. 
The subgroups of $\widetilde D_n$ are $A_{q,0}$, 
$A_{q,0}\cup B_{q,r}$, and the left cosets are 
$A_{q,r}$, $B_{q,r}$, $A_{q,r}\cup B_{q,s}$ 
($q$ a positive divisor of $2n$, $r,s\in\{0,1,\ldots,q-1\}$). 

Take now $m+1\geq 5$ cosets $g_iG_i$ $(i=0,1,\ldots,m)$
in $\widetilde D_n$, and assume that 
any $4$ of them have a common element. 
Partition each given coset as $g_iG_i=A^i\cup B^i$, where $A^i$ is either empty or of the form $A_{q,r}$, and $B^i$ is either empty or of the form $B_{q,r}$. 

Case I. $A^i\cap A^j$ is non-empty for all $i<j$: 
Then by the proof of Proposition~\ref{prop:cyclic} we have that 
the intersection of all the $A^i$ is non-empty, and we are done. 

Case II. There two among the $A^i$, say $A^0$ and $A^1$ such that 
$A^0\cap A^1=\emptyset$. Then $g_0G_0\cap g_1G_1=B^0\cap B^1$, hence for any 
$2\leq i<j\leq m$ we have that 
$g_0G_0\cap g_1G_1\cap g_iG_i\cap g_jG_j=B^0\cap B^1\cap B^i\cap B^j$ 
is non-empty. 
Write $q_j$ for the positive divisor of $2n$ and 
$r_j\in \{0,1,\ldots,q_j-1\}$ with 
$B^j=B_{q_j,r_j}$, $j=0,\ldots,m$. 
Consider the following system of congruences: 
\[x\equiv r_0\mbox{ mod }q_0;\quad x\equiv r_1\mbox{ mod }q_1;\quad 
x\equiv r_2\mbox{ mod }q_2;\quad\ldots;\quad x\equiv r_m\mbox{ mod }q_m.\]
Our assumption implies that for any $2\leq i<j\leq m$, the subsystem 
consisting of the congruences with indices $0,1,i,j$
has a solution. Consequently, 
any pair of the above $m+1$ congruences has a simultaneous solution. Therefore 
the whole system has a solution by the Chinese Remainder Theorem. 
In other words, the $m+1$ cosets have a non-empty intersection. 
\end{proof}

\begin{corollary}\label{cor:helly}
The Helly dimension of a finite subgroup of $\sltwo$ is bounded by $6$. 
\end{corollary} 

\begin{proof} 
An inspection of the corresponding subgroup lattices shows that 
$\lambda(A_4)=3$, $\lambda(S_4)=4$, and $\lambda(A_5)=4$. 
It is easy to see that if $G\cong \widetilde G/N$, where $N$ is a two-element normal subgroup in the finite group $\widetilde G$, then $\lambda(\widetilde G)=\lambda(G)+1$. Hence we have $\lambda(\widetilde A_4)=4$, 
$\lambda(\widetilde S_4)=5$, and $\lambda(\widetilde A_5)=5$. 
By Lemma~\ref{lemma:inequalities} we get the desired bound for the Helly dimension of these groups. The two infinite series of finite subgroups of 
$\sltwo$ are dealt with in Propositions~\ref{prop:cyclic} and 
\ref{prop:dihedral}. 
\end{proof}


\section{Binary forms}\label{sec:binary}

Throughout this section our group $G$ is the special linear group 
$\sltwo$, $V=\pol_d(\mc^2)$ is the space of homogeneous binary forms of degree $d$ with complex coefficients, endowed with the natural action 
of $\sltwo$, whereas $W$ is an arbitrary finite dimensional rational $\sltwo$--module. Recall that 
$\pol_d(\mc^2)$, $d=0,1,2,\ldots$ is a complete list of representatives of the isomorphism classes of irreducible rational 
$G$--modules. 
The study of the invariants of several binary forms is the most classical 
topic in invariant theory. 
Our aim here is to investigate the number of type $V$ 
variables needed in a typical system of separating invariants in this case. 
It turns out that there is a uniform upper bound on $\sigma(G,W,V)$ 
(independent of $W$ and $V$). 

\begin{theorem}\label{thm:binary}
We have $\sigma(\sltwo,W,V)\leq 8$ for any finite dimensional $\sltwo$--module 
$W$ and irreducible $\sltwo$--module $V$; that is, 
if $S$ is a separating subset in $\mc[W\oplus V^8]^{\sltwo}$ 
(eight copies of $V$), then $S$ is a typical system of separating invariants in type $V$ variables relative to $W$. 
\end{theorem} 

{\bf Remark.} An arbitrary finite dimensional $\sltwo$--module is isomorphic 
to $V_1^{m_1}\oplus\cdots\oplus V_r^{m_r}$ 
(where $V_1,\ldots,V_r$ are pairwise non-isomorphic irreducible 
$\sltwo$--modules). Theorem~\ref{thm:binary} implies that 
$\mc[\bigoplus_{i=1}^rV_i^{m_i}]^{\sltwo}$ 
contains a finite separating set whose elements depend only on $\leq 8$ 
variables of type $V_i$ for each $i=1,\ldots,r$. 

As we mentioned already in Section ~\ref{sec:reductive}, two points in a 
$G$--module can be separated by a polynomial invariant if and only if 
the unique closed orbits in their 
orbit closures are different. So Theorem~\ref{thm:binary} is equivalent to the 
following. 

\begin{proposition}\label{prop:3}
Let $v=(w,v_1,\ldots,v_m)$ and $v'=(w',v_1',\ldots,v_m')$ be points in 
$W\oplus V^m$ whose $\sltwo$--orbits are closed. 
Assume that for all $0\leq r\leq 8$, $1\leq i_1<\cdots<i_r\leq m$, 
the  closures of the orbits of $(w,v_{i_1},\ldots,v_{i_r})$ 
and $(w',v_{i_1}',\ldots,v_{i_r}')$ in $W\oplus V^r$ 
have a non-empty intersection. 
Then $v$ and $v'$ belong to the same orbit. 
\end{proposition} 

In order to prove this statement we shall extensively use 
the generalization of the 
Hilbert-Mumford Criterion due to Birkes-Richardson 
(see \cite{b} or Theorem 6.9 in \cite{pv}). 

{\it BRHM Criterion}: Assume that the $G$--orbit of $u_1\in W$ is closed. 
Then $u_1$ belongs to the Zariski closure of the $G$--orbit of $u_2$ 
if and only if there exists a {\it one-parameter subgroup} 
(shortly $\onepsg$) $\rho:\mc^{\times}\to G$ 
(a homomorphism of algebraic groups) such that $\limzr u_2$ exists and belongs to the $G$--orbit of $u_1$. 

Any non-trivial $\onepsg$ in $\sltwo$ is of the following form: 
there is an element $g\in\sltwo$ and a natural number $n$ such that 
\begin{equation}\label{eq:1psg}
\rho(z)=g\left(\begin{array}{cc}z^{-n} & 0   \\0 & z^n 
\end{array}\right)
g^{-1}.
\end{equation}
We say that a $\onepsg$ 
{\it acts on a linear form $l\in \pol_1(\mc^2)$ by a positive character} 
if the line $\mc l$ is stabilized by $\rho$ and $\rho(z)l=z^nl$ 
for a positive integer $n$. Note that if $\rho$ is of the form 
(\ref{eq:1psg}), then the only lines fixed by $\rho$ are spanned by  
$gx$ and $gy$ (where $x$ and $y$ stand for the usual coordinate functions on $\mc^2$), 
and $\rho$ acts by a positive character on $gx$ (but not on $gy$). 
An orbit is called {\it maximal} if it is not contained in the closure 
of another orbit. 
We say that the linear form $l$ is {\it a root of multiplicity} $n$ 
of the binary form $v$ if $l^n$ divides $v$, and $l^{n+1}$ does not divide $v$. 
From the BRHM Criterion 
and (\ref{eq:1psg}) one derives the following known facts. 

\begin{proposition}\label{prop:root} 
\begin{itemize}\item[(i)] Let $\rho$ be a $\onepsg$ as in (\ref{eq:1psg}), 
and $v$ a homogeneous binary form of degree $d$. 
Then $\limzr v$ exists if and only if $gx$ is a root of multiplicity 
$\geq d/2$ of $v$. If the multiplicity is greater than $d/2$, then 
$\limzr v=0$, and if $d$ is even and the multiplicity is $d/2$, then 
$\limzr v$ is a non-zero scalar multiple of $(gx)^{d/2}(gy)^{d/2}$. 

\item[(ii)] If $\limzr v$ exists for some non-trivial $\onepsg$, then 
it does not belong to the $\sltwo$--orbit of $v$, unless $v=0$ or 
$d$ is even and $v=(l_1l_2)^{d/2}$ for some independent linear forms $l_1,l_2$. 
\end{itemize}
\end{proposition}

\begin{proposition}\label{prop:maximal} 
The $\sltwo$--orbit of $v=(w,v_1,\ldots,v_m)\in W\oplus V^m$ is closed 
and maximal if the components $v_1,\ldots,v_m$ have no common root 
of multiplicity $\geq d/2$ (in all $v_i$). 
In particular, the orbit of $v$ is closed and maximal if 
there is a type $V$ component that has no root of muliplicity 
$\geq d/2$, or if  
there are two type $V$ components 
that have no common root of multiplicity $\geq d/2$. 
\end{proposition}

\begin{proof} Assume that the orbit of $v$ is not closed. 
By the BRHM Criterion, there is a non-trivial $\onepsg$ $\rho$ 
such that $\limzr v$ exists. Then $\limzr v_i$ exists for all $i$. 
We may assume that $\rho$ is of the form 
(\ref{eq:1psg}). Then $gx$ is a root of all $v_i$ with multiplicity $\geq d/2$ 
by Proposition~\ref{prop:root}. 
Next, assume that the orbit of $v$ is not maximal. 
Then there is some 
$u\in W\oplus V^m$ and a non-trivial $\onepsg$ $\rho$ (of the form 
(\ref{eq:1psg})) with $v=\limzr u$, implying 
that $v_i=\limzr u_i$ for all $i$. If $d$ is odd, then all the $v_i$ are 
zero, and if $d$ is even, then all $v_i\in \mc (gx)^{d/2}(gy)^{d/2}$ 
by Proposition~\ref{prop:root}, hence say $gx$ is a common root of the $v_i$ 
with multiplicity $d/2$. 
\end{proof}  

Similar arguments yield the following consequence of Proposition~\ref{prop:root}. 

\begin{proposition}\label{prop:notinclosure}
Assume that the orbit of $v=(w,v_1,\ldots,v_m)\in W\oplus V^m$ is closed. 
Then it is maximal if either 
of the following conditions holds: 
\begin{itemize}
\item[(i)] $v$ has a type $V$ component which is non-zero, and can not be written as $(l_1l_2)^{d/2}$ where $d$ is even and $l_1,l_2$ are independent linear forms;   
\item[(ii)] $v$ has two linearly independent components of type $V$. 
\end{itemize}
\end{proposition}

{\bf Remark.} A description of the orbit cosures 
is given in \cite{kh} for irreducible $\sltwo$--modules 
and in  \cite{p} for arbitrary $\sltwo$--modules. 
One could deduce Propositions~\ref{prop:maximal} and \ref{prop:notinclosure} 
from the results of \cite{p}.

We shall need also the following elementary fact on the irreducible representations of $\sltwo$. 

\begin{lemma}\label{lemma:stabilizer}
The pointwise stabilizer of a linear subspace of dimension $\geq 2$ 
in the irreducible $\sltwo$--module $V$ is finite. 
\end{lemma} 

\begin{proof} Let $H$ be the pointwise stabilizer in $\sltwo$ of a subspace of $V$ of dimension $\geq 2$. 
Any element of $H$ is conjugate in $\sltwo$ to one of 
$\left(\begin{array}{cc}1 & 1 \\0 & 1\end{array}\right)$ 
or 
$\left(\begin{array}{cc}-1 & 1 \\0 & -1\end{array}\right)$ 
or  
$\left(\begin{array}{cc}z & 0 \\0 & z^{-1}\end{array}\right)$ 
$(z\in\mc^{\times})$. 
Direct computation shows that the fixed point subspaces of the unipotent elements are $1$--dimensional. 
Therefore any element of $H$ is semisimple. It follows that the connected  component of the identity in $H$ is contained in a maximal torus $T$ of 
$\sltwo$. Now $T$ is $1$--dimensional (as an algebraic variety), and 
the fixed point subspace of $T$ is $1$--dimensional if $d$ is even and $0$--dimensional if $d$ is odd. Consequently, the connected componenent of the identity in $H$ is $0$--dimensional, hence is the trivial subgroup 
$\{1\}$, showing that $H$ is finite. 
\end{proof}

\bigskip
\begin{proofofprop}~\ref{prop:3}. If all $v_i$, $v_i'$ are zero, then the assertion holds trivially. From now on we assume that this is not the case. 

{\it Case I.} There is a type $V$ component of $v$, say $v_1$, that has no root of multiplicity $\geq d/2$. 

We claim that for all $0\leq r\leq 7$ and $2\leq i_1<\cdots<i_r\leq m$, 
the $G$--orbits of $(w,v_1,v_{i_1},\ldots,,v_{i_r})$ and 
$(w',v_1',v_{i_1}',\ldots,v_{i_r}')$ coincide. 
Indeed, the closures of these orbits intersect by assumption. On the other hand, the first $G$--orbit is closed and maximal by Proposition~\ref{prop:maximal}. 
These two facts clearly imply that the orbits coincide. 

{\it Case I/a.} $v_i=\alpha_i v_1$ with $\alpha_i\in\mc$ for 
$i=2,\ldots,m$. 

Since $(w,v_1,v_i)$ and $(w',v_1',v_i')$ belong to the same $G$--orbit, we have also $v_i'=\alpha_i v_1'$ for $i=2,\ldots,m$. Take $g\in G$ with $g(w,v_1)=
(w',v_1')$, such a $g$ exists, and 
$gv_i=g(\alpha_iv_1)=\alpha_i gv_1= \alpha_iv_1'=v_i'$ for $i=2,\ldots,m$, 
so $gv=v'$. 

{\it Case I/b.} There is an $i\geq 2$, say $i=2$ such that $v_1$ and $v_2$ 
are linearly independent. 

Our assumptions imply that $u=(w,v_1,v_2)$ and $u'=(w',v_1',v_2')$ belong to the same 
$G$--orbit, so replacing $v$ by an appropriate element in its orbit 
we may assume that $u=u'$. By the assumptions, for all $1\leq r\leq 6$ and 
$3\leq i_1<\cdots<i_r\leq m$, the elements $(u,v_{i_1},\ldots,v_{i_r})$ 
and $(u,v_{i_1}',\ldots,v_{i_r}')$ belong to the same $G$--orbit, 
hence $(v_{i_1},\ldots,v_{i_r})$ and $(v_{i_1}',\ldots,v_{i_r}')$ 
belong to the same orbit with respect to the stabilizer $H$ of $u$. 
Now $H$ is a finite subgroup of $\sltwo$ by 
Lemma~\ref{lemma:stabilizer}, 
so we can conclude from Corollary~\ref{cor:helly} that 
$(v_3,\ldots,v_m)$ and $(v_3',\ldots,v_m')$ belong to the same $H$--orbit, 
or in other words, that $v$ and $v'$ belong to the same $G$--orbit. 

{\it Case II.} All components of $v$ of type $V$ have a root of multiplicity 
$\geq d/2$, and there are two components, say $v_1$ and $v_2$ that have no common root of multiplicity $\geq d/2$ (both in $v_1$ and $v_2$). 

The orbit of $(w,v_1,v_2,v_{i_1},\ldots,v_{i_r})$ is closed and maximal 
by Proposition~\ref{prop:maximal}, hence coincides with the orbit of 
$(w',v_1',v_2',v_{i_1}',\ldots,v_{i_r}')$ for all $0\leq r\leq 6$ and 
$3\leq i_1<\cdots <i_r\leq m$. 
Moreover, $v_1$ and $v_2$ are linearly independent. One concludes that 
the orbits of $v$ and $v'$ are the same just like in Case I/b. 

{\it Case III.} There is a linear form $l$ such that $l^{e+1}$ divides all of $v_1,\ldots,v_m$, where $d=2e$ or $d=2e+1$, and not all $v_i$ are zero. 

Note that $l$ is the only root of multiplicity $\geq d/2$ for each non-zero 
component $v_i$. 
Without loss of generality we may assume that $v_1$ is non-zero. 
By Proposition~\ref{prop:root}, 
$\limzr (v_1,\ldots,v_m)$ exists if and only if $\rho$ acts by a positive character on $l$, and in this case $\limzr(v_1,\ldots,v_m)$ does not belong to the orbit of $(v_1,\ldots,v_m)$. Hence in this case $\limzr w$ does not exist 
(recall that the orbit of $v$ is closed by assumption). Consequently, 
$\limzr (w,v_1)$ does not exist for any non-trivial 
$\onepsg$, implying that 
the orbit of $(w,v_1,v_{i_1},\ldots,v_{i_r})$ is closed for all 
$r\geq 0$ and $2\leq i_1<\cdots<i_r\leq m$. Moreover, these orbits are maximal 
by Proposition~\ref{prop:notinclosure}. 
Thus the assumptions of the Proposition imply that 
$(w,v_1,v_{i_1},\ldots,v_{i_r})$ and 
$(w',v_1',v_{i_1}',\ldots,v_{i_r}')$ belong to the same orbit 
provided that $r\leq 7$. One finishes in the same way as in Case I. 

Since the roles of $v$ and $v'$ are symmetric, one can deal with 
Cases I',II',III' (obtained by replacing $v$ by $v'$) in the same way as 
with Cases I,II,III. 

Note that the proof is already complete for odd $d$. 
There are however further cases to consider when $d$ is even, which we 
assume from now on. 
From now on we automatically assume that we are not in 
the cases covered so far. In particular, this implies that 
any two type $V$ components of $v$ have a common root of multiplicity 
$\geq d/2$, and the same holds for $v'$. 

{\it Case IV.} There is a non-zero  type $V$ component of $v$ or $v'$, say $v_1$, 
which can not be written as $(l_1l_2)^{d/2}$ with two independent linear forms 
$l_1$ and $l_2$. 

There is a unique (up to scalar) root $l$ of multiplicity $\geq d/2$ of $v_1$, 
and $l^{d/2}$ necessarily divides all the $v_i$. 
Since the orbit of $v$ is closed, just like in Case III, we conclude from  Propositions~\ref{prop:root} and \ref{prop:notinclosure} that 
the orbit of 
$(w,v_1,v_{i_1},\ldots,v_{i_r})$ is closed and maximal for all $0\leq r$ and 
$2\leq i_1<\cdots<i_r\leq m$. 
One finishes in the same way as in Case I (or III). 

From now on we assume also that we are not in Case IV. 
This means in particular that each non-zero type $V$ component of 
$v$ and $v'$ can be written as the product of the $(d/2)$th powers of 
two independent linear forms. 

{\it Case V.} There is a non-zero linear form $l$ and linear forms $l_1,\ldots,l_m$ such that 
$v_i=l^{d/2}l_i^{d/2}$ for $i=1,\ldots,m$, and there are two among the $l_i$, 
say $l_1$ and $l_2$ that are linearly independent. 

For any $\onepsg$ $\rho$ acting on $l$ by a positive character
$\limzr (v_1,\ldots,v_m)$ exists, furthermore,  
$\limzr (v_1,v_2)$ does not belong to the orbit of $(v_1,v_2)$. 
Since the orbit of $v$ is closed, this implies $\limzr w$ does not exist 
for such $\rho$. 
Taking into account that $l$ is the only common root of $v_1$ and $v_2$, this implies that 
the orbit of $(w,v_1,v_2,v_{i_1},\ldots,v_{i_r})$ is closed 
for any $r\geq 0$ and $3\leq i_1<\cdots <i_r\leq m$. These orbits are maximal by 
Proposition~\ref{prop:notinclosure}. 
One can finish as in Case I/b. 

{\it Case V'.} Replace $v$ by $v'$ in Case V. 

{\it Case VI.} $v_i=\alpha_i (l_1l_2)^{d/2}$, 
$v_i'=\beta_i (b_1b_2)^{d/2}$, $i=1,\ldots,m$, where 
$\alpha_i,\beta_i\in \mc$, the linear forms $l_1$, $l_2$ are independent, 
and the linear forms $b_1$, $b_2$ are independent. 

There is a non-zero scalar among the $\alpha_i$,$\beta_j$, say 
$\alpha_1\neq 0$. 
Then the $G$--orbit of $(w,v_1)$ in $W\oplus V$ 
is closed by Lemma~\ref{lemma:closed}. 
If $v_1'=0$, then $(w,v_1)$ is not contained in the closure 
of the orbit of $(w',v_1')$, contradicting our assumption in the Proposition. 
So $\beta_1\neq 0$, and the orbit of $(w',v_1')$ is closed. 
Moreover, the orbits of $(w,v_1,v_{i_1},\ldots,v_{i_r})$ and 
$(w',v_1',v_{i_1}',\ldots,v_{i_r}')$ are also closed by Lemma~\ref{lemma:closed} 
for all $r\geq 0$ and $2\leq i_1<\cdots<i_r\leq m$, hence they coincide by our assumptions if $r\leq 7$. One finishes in the same way as in 
Case I/a. 

For the rest of the proof we assume additionally that we are not in 
Cases V,V',VI. So $v_1,\ldots,v_m$ have no common root of multiplicity 
$\geq d/2$. On the other hand, any pair $v_i,v_j$ have a common root of multiplicity $\geq d/2$. Moreover, the same holds for $v'$. 
This implies that there are three pairwise independent linear forms $l_1,l_2,l_3$ such that a non-zero scalar multiple of 
each of $(l_1l_2)^{d/2}$, $(l_1l_3)^{d/2}$, $(l_2l_3)^{d/2}$ 
occurs among the components, say 
$v_1=(l_1l_2)^{d/2}$, $v_2=(l_1l_3)^{d/2}$, $v_3=(l_2l_3)^{d/2}$, and all 
the $v_i$ belong to $\mc v_1\cup\mc v_2\cup\mc v_3$. 
By Proposition~\ref{prop:maximal}, 
the orbit of 
$(w,v_1,v_2,v_3)$ is closed and maximal, 
hence by assumption coincides with the orbit of 
$(w',v_1',v_2',v_3')$, and so we may assume that they are equal 
(replace $v$ by an appropriate element in its orbit). 
It follows that all type $V$ components of 
$v'$ belong to $\mc v_1\cup\mc v_2\cup\mc v_3$. Take a non-zero component $v_i$ of $v$, where $i\geq 4$. 
Then there is a nonzero $\alpha_i\in\mc$ and $f(i)\in\{1,2,3\}$ 
such that $v_i=\alpha_iv_{f(i)}$. Assume for example that $f(i)=3$. 
Then the same argument as above shows that $(w,v_1,v_2,v_i)$ and 
$(w',v_1',v_2',v_i')=(w,v_1,v_2,v_i')$ belong to the same orbit. 
Therefore $v_i'\notin \mc v_1\cup\mc v_2$, so 
$v_i'=\beta_i v_3$ with a non-zero $\beta_i\in \mc$. 
By assumption, the orbit closures of $(v_3,v_i)$ and $(v_3',v_i')$ have a non-empty intersection, implying that $\alpha_i=\beta_i$, hence $v_i=v_i'$. 
By symmetry in $v$ and $v'$, we have that $v_i$ is non-zero if and only if 
$v_i'$ is non-zero. So a repeated use of the above argument implies that $v_i=v_i'$ for all $i\geq 4$, hence $v=v'$. 
\end{proofofprop}

We take a digression and state and prove the analogues of Theorem~\ref{thm:binary} for ordinary polarization. 

\begin{theorem}\label{thm:polarization} 
If $S$ is a separating set of invariants in $\mc[W\oplus V^3]^{\sltwo}$ 
(three copies of $V$), then the polarizations 
with respect to the type $V$ variables 
of the elements of $S$ 
constitute a separating set in $\mc[W\oplus V^m]^{\sltwo}$ for an arbitrary 
$m$.  
Moreover, it is sufficient to consider polarizations of the following special form: when $m\geq 3$, 
take the coefficients of $f(w,v_{i_1},v_{i_2},\sum_{j=3}^{{\mathrm{min}}\{m,8\}}t^{j-3}v_{i_j})$ 
(viewed as a polynomial in $t$), as $f$ ranges over $S$, and 
$i_1,i_2,i_3,\ldots$ are different elements of $\{1,\ldots,m\}$; 
when $m<3$, take $f(w,v_1,\ldots,v_m,0,\ldots,0)$ as $f$ ranges over $S$.  
\end{theorem} 

\begin{proof} The proof is almost the same as the proof of Theorem~\ref{thm:binary}. Write $q$ for ${\mathrm{min}}\{m,8\}$. We treat the same Cases as in that proof. The only difference is 
that when we reduced a given Case to the study of orbits of a finite group $H$ 
(the stabilizer of $(w,v_{i_1},v_{i_2})$, then we need to point out that if 
$(v_{i_3},\ldots,v_{i_q})$ 
and 
$(v_{i_3}',\ldots,v_{i_q}')$ 
belong to different $H$--orbits, then there is a $t\in\mc$ such that $\sum_{j=3}^qt^{j-3}v_{i_j}$ belongs to a different $H$--orbit than 
$\sum_{j=3}^qt^{j-3}v_{i_j}'$; this folllows from the results on "cheap polarization" in \cite{dkw}. Hence the closed maximal orbit of 
$(w,v_{i_1},v_{i_2},\sum_{j=3}^qt^{j-3}v_{i_j})$ is not contained in the closure of the orbit of $(w',v_{i_1}',v_{i_2}',\sum_{j=3}^qt^{j-3}v_{i_j}')$, implying 
that 
$f(w,v_{i_1},v_{i_2},\sum_{j=3}^qt^{j-3}v_{i_j})\neq 
f(w',v_{i_1}',v_{i_2}',\sum_{j=3}^qt^{j-3}v_{i_j}')$ 
for an appropriate $f\in S$. 
\end{proof} 

{\bf Remark.} The bound $3$ on the number of type $V$ variables is sharp, 
as it can be seen already when $\dim(V)=3$.


\section{Rational invariants}\label{sec:rational} 

Throughout this section $k$ is algebraically closed. 
Our results on typical systems of separating invariants have a version for rational invariants as well. 
Given an irreducible $G$--variety $X$, 
denote by $k(X)$ the field of rational functions on $X$ and by 
$k(X)^G$ its subfield of invariants. 
Following \cite{pv} we say that 
$R\subset k(X)^G$ {\it separates orbits in general position} 
if there is a non-empty open subset $X_0$ of $X$ such that if 
$x,x'\in X_0$ belong to different orbits, then there is a rational invariant 
$f\in R$ defined both at $x$ and $x'$ such that $f(x)\neq f(x')$. 
It is well known that a finite set $R$ of rational invariants 
separates orbits in general position 
if and only if $k(X)^G$ is a purely inseparable extension of its subfield generated by $k$ and $R$; see \cite{pv} for the characteristic zero case, 
and combine it with (18.2) from \cite{bo} to get the general statement. 
A theorem of Rosenlicht \cite{r} asserts that a finite set $R$ of rational 
invariants separating orbits in general position always exists. 
(Note that if a finite set $R$ of rational invariants separates orbits in general position, then there is a non-empty $G$--stable open subset $U$ in $X$ such that the elements of $R$ are all defined on $U$, and separate 
inequivalent points in $U$, thus a geometric quotient $U/G$ exists.) 

Write $X^m$ for the product $X\times\cdots\times X$ of $m$ copies of $X$. 

\begin{theorem}\label{thm:rational} 
Let $X$ be an irreducible $G$--subvariety of $V=k^n$ (where $G$ acts linearly on $V$), let $Y$ be an arbitrary 
irreducible $G$--variety, and assume that $R\subset k(Y\times X^{n+1})^G$ separates orbits in general position in $Y\times X^{n+1}$.  
Then for any $m\geq n+1$ we have that 
\[R^{(m)}=\{f^{(i_1,\ldots i_{n+1})}\mid 1\leq i_1<\cdots<i_{n+1}\leq m\}\] 
separates orbits in general position in $Y\times X^m$. 
In particular, if $\char(k)=0$, then $R^{(m)}$ generates the field of rational 
invariants on $Y\times X^m$. 
\end{theorem}

\begin{proof} 
Let $U$ be a non-empty open subset of $Y\times X^{n+1}$ such that 
inequivalent points in $U$ can be separated by an element of $R$. 
The complement of $U$ in $Y\times X^{n+1}$ is a proper closed 
subset $Z$ in $Y\times X^{n+1}$. 
For $1\leq i_1<\cdots <i_{n+1}\leq m$, denote by 
$\pi_{(i_1,\ldots,i_{n+1})}$ the projection morphism 
$Y\times X^m\to Y\times X^{n+1}$, 
$(y,x_1,\ldots,x_m)\mapsto (y,x_{i_1},\ldots,x_{i_{n+1}})$. 
Then 
\[Z^{(m)}=\bigcup_{1\leq i_1<\cdots<i_{n+1}\leq m}
\pi_{(i_1,\ldots,i_{n+1})}^{-1}(Z)\] 
is a proper closed subset of $Y\times X^m$; 
denote its non-empty complement by $U^{(m)}$. 
(If $U$ was $G$--stable, then $U^{(m)}$ is $G$--stable, since the 
projection maps are $G$--equivariant.)  
Note that if $x\in U^{(m)}$, then any projection 
$\pi_{(i_1,\ldots,i_{n+1})}(x)$ belongs to $U$. 
If for $x,x'\in U^{(m)}$ there is an $(i_1,\ldots,i_{n+1})$ such that 
$y=\pi_{(i_1,\ldots,i_{n+1})}(x)$ and 
$y'=\pi_{(i_1,\ldots,i_{n+1})}(x')$ belong to different $G$--orbits, 
then $y$ and $y'$ can be separated by some $f\in R$, and so 
$f^{(i_1,\ldots,i_{n+1})}$ separate $x$ and $x'$. 

Assume now that all projections $\pi_{(i_1,\ldots,i_{n+1})}$ map 
$x=(y,x_1,\ldots,x_m)\in U^{(m)}$ and 
$x'=(y',x_1',\ldots,x_m')$ 
into the same $G$--orbit in $Y\times X^{n+1}$. 
After a possible renumbering of the indices, we may assume that 
$x_1,\ldots,x_n$ span the $k$--subspace $\span_k\{x_1,\ldots,x_m\}$ 
of $V$. So there are coefficients $\alpha_{ij}$ with 
$x_{n+i}=\sum_{j=1}^n\alpha_{ij} x_j$ for $i=1,\ldots,m-n$. 
Since $(y,x_1,\ldots,x_n,x_{n+i})$ and $(y',x_1',\ldots,x_n',x_{n+i}')$ 
belong to the same $G$--orbit of $Y\times X^{n+1}$ for all $i$ 
by our assumption, 
we have that 
$x_{n+i}'=\sum_{j=1}^n\alpha_{ij}x_j'$ holds for $i=1,\ldots,m-n$ 
(recall that the action of $G$ on $X$ is induced by a linear action on $V$). 
It follows that for an arbitrary element $g\in G$ with 
$g\cdot(y,x_1,\ldots,x_n)=(y',x_1',\ldots,x_n')$ we have that $g\cdot x=x'$. 
Such a $g$ exists by our assumption, hence $x$ and $x'$ belong to the 
same $G$--orbit of $Y\times X^m$. 
\end{proof}


\section{Concluding remarks}\label{sec:conclude} 

The present paper is motivated by Weyl's Theorem \cite{w}, Theorem 2.5A, 
asserting that 
if $k$ is of characteristic zero, then $k[W\oplus V^m]^G$ 
is generated by the polarizations (with respect to the type $V$ variables) 
of the elements of a generating system 
of $k[W\oplus V^n]^G$, where $n=\dim(V)$. 
It is well known that this does not hold when $\char(k)$ is positive (see for example \cite{do} for a class of examples).  
However, in a recent preprint, J. Draisma, G. Kemper and D. Wehlau \cite{dkw} 
proved that polarizing the elements of a separating set in 
$k[W\oplus V^n]^G$, $n=\dim(V)$, one always gets a separating 
set in $k[W\oplus V^m]^G$, regardless of the characteristic of the base field. 

Our statements Theorem~\ref{thm:2n} and ~\ref{thm:n+1} 
are weaker if one is interested in degree bounds for separating invariants. 
On the other hand, they are "really" about separating invariants, in the sense that they have no analogues for generating invariants, even in characteristic zero. 
From the point of view of polarizations, our results can be viewed as an 
addendum to \cite{dkw} making possible some computational simplifications in the construction 
of a separating set in $k[W\oplus V^m]^G$ for large $m$. 
Namely, one takes a separating set $S_0$ in $k[W\oplus V^n]^G$ with $n=\dim(V)$. 
Then it is sufficient to take those polarizations of the elements of $S_0$ 
that depend only on $\leq 2n$ (resp. $\leq n+1$ for reductive $G$) type $V$ variables to get a (multihomogeneous) 
separating set in $k[W\oplus V^m]^G$.   
This is an essential simplification from the point of view of explicit computations. 

Finally, this trivial version of "polarization" applies also for rational invariants, where the general polarization process does not make sense. 

\begin{center}{\bf Acknowledgement}\end{center} 

I thank P. P. P\'alfy for remarks that led to an improved presentation of 
Section~\ref{sec:finite}. 
This research was partially supported by OTKA No. T046378, the Bolyai Fellowship 
and Leverhulme Research Exchange Grant F/00158/X.


 \end{document}